\documentclass[12pt]{article}
\usepackage{enumerate}
\usepackage{amsopn,amsfonts,amsmath,amssymb}
\usepackage{float}
\usepackage[american]{babel}
\usepackage{color, colortbl}
\usepackage{float}
\usepackage{dsfont}
\usepackage{color, colortbl}
 \usepackage[utf8]{inputenc}
 \usepackage{natbib}

\usepackage{setspace}
\allowdisplaybreaks
\usepackage{caption}

\numberwithin{equation}{section}

\def\R{{\mathbb R}}

\def\N{{\mathbb N}}

\def\rank{{\rm rank}}

\def\P{\mathbb{P}}

\def\op{\mathrm{op}}
\usepackage[colorlinks=true,
            linkcolor=red,
            urlcolor=blue,
						citecolor=black]{hyperref}

\newtheorem{theorem}{Theorem}[section]
\newtheorem{lemma}{Lemma}[section]

\newtheorem{corollary}{Corollary}[section]

\newtheorem{assumption}{Assumption}[section]

\DeclareMathOperator*{\argmin}{argmin}
\DeclareMathOperator*{\argmax}{argmax}
\usepackage[left=2.5cm, right=2.5cm, top=2.5cm]{geometry}

\numberwithin{equation}{section}

\begin{document}

\begin{titlepage}

\title{\Large \bf Factor and factor loading augmented estimators for panel regression}

\author{
{\large Jad B\textsc{eyhum}}\footnote{
ORSTAT, KU Leuven, \href{mailto:jad.beyhum@kuleuven.be}{jad.beyhum@kuleuven.be}. This work was undertaken at the Toulouse School of Economics, Université Toulouse Capitole.}
 \and
{\large Eric G\textsc{autier}}\thanks{Toulouse School of Economics, Université Toulouse Capitole, \href{eric.gautier@tse-fr.eu}{eric.gautier@tse-fr.eu}.\newline
 Financial support from the European Research Council (2014-2019 / ERC grant agreement No.\ 337665) is gratefully acknowledged. The authors thank Domenico Giannone, Jihyun Kim, Pascal Lavergne, Thierry Magnac and Nour Meddahi for helpful comments and ideas.}}

\date{}

\maketitle
\begin{abstract}
This paper considers linear panel data models where the dependence of the regressors and the unobservables is modelled through a factor structure. The asymptotic setting is such that the number of time periods and the sample size both go to infinity.  Non-strong factors are allowed and the number of factors can grow to infinity with the sample size. We study a class of two-step estimators of the regression coefficients. In the first step, factors and factor loadings are estimated. Then, the second step corresponds to the panel regression of the outcome on the regressors and the estimates of the factors and the factor loadings from the first step. Different methods can be used in the first step while the second step is unique. We derive sufficient conditions on the first-step estimator and the data generating process under which the two-step estimator is asymptotically normal. Assumptions under which using an approach based on principal components analysis in the first step yields an asymptotically normal estimator are also given. The two-step procedure exhibits good finite sample properties in simulations.
\end{abstract}
\vspace{0in}
\end{titlepage}

\newpage 
\section{Introduction}
%Panel data allows to estimate regression models where the error term exhibits cross-section and and serial correlations. 
This paper considers inference on $\beta\in \mathbb{R}^K$ in the following model: \begin{equation}
Y_{it} = \sum_{k=1}^K\beta_kX_{kit} +\sum_{j=1}^{r_N}\lambda_{ij}f_{tj}\delta_{j}+E_{it},\label{model 1},
\end{equation}
where %$i\in\{1,...,N\}$ indices the individuals and $t\in\{1,...,T\}$ the time periods, 
the data consists of the outcome $Y_{it}$ and the regressors $X_{kit}$ for all $k=1,\dots,K$, $i=1,\dots,N$ and $t=1\dots,T$. The random vectors $\lambda_i$ and $f_t$ in $\mathbb{R}^{r_N}$ are factor loadings and factors, $\delta$ is a nonrandom vector in $\mathbb{R}^{r_N}$, $r_N$ is the number of factors.

This is a panel data model with interactive fixed effects (see \citet{pesaran2015time}). It allows for flexible cross-section and serial correlation thanks to the factor structure in the regression error. Several techniques have been developed to estimate this model.  \citet{bai2009panel} proposes to estimate jointly the regression coefficient and the factors and factor loadings. \citet{moon2018nuclear} and \citet{beyhum2019square} study a nuclear-norm penalized estimator.  In contrast, the CCE estimator of \citet{pesaran2006estimation} and the factor-augmented regression estimator studied in \citet{kapetanios2005alternative} and \citet{greenaway2012asymptotic} model the dependence between the regressors and the unobservables $\sum_{j=1}^{r_N}\lambda_{ij}f_{tj}\delta_{j}+E_{it}$. They assume that, for $k\in\{1,\dots, K\}$, there exists $\lambda_{k1},\dots,\lambda_{kN}$ which are random vectors in $\mathbb{R}^{r_N}$ and mean-zero errors $E_1,\dots,E_K$ which are $N\times T$ random matrices such that
$X_{kit} = \sum_{r=1}^{r_N}\lambda_{kir}f_{tr}+E_{kit}$ for $k\in\{1,\dots,K\}.$ This means that the regressors have a factor structure with the same factors as the error term but possibly different factor loadings. 

In the papers of \citet{pesaran2006estimation}, \citet{bai2009panel} and \citet{greenaway2012asymptotic}, a strong factor assumption is imposed. It means that the ratio of the singular values of $\Gamma$ and $\sqrt{NT}$ has a finite deterministic limit as $N,T\to \infty$, where $\Gamma_{it}=\sum_{r=1}^{r_N}\lambda_{ij}f_{tj}$. It holds if $(1/N)\sum_{i=1}^N\lambda_i\lambda_i^\top$ $(1/T)\sum_{t=1}^Tf_tf_t^\top$ has a finite deterministic limit in probability. The number of factors is also assumed to be fixed with the sample size. It is worth noting that some papers have sought to relax these assumptions in the context of the CCE (\citet{chudik2011weak}) and factor augmented (\citet{reese2018estimation}) estimators.

This paper proposes instead to model the dependence of the regressors with both the factors and the factor loadings by assuming that there exists $\delta_k\in\mathbb{R}^{r_N}$ for $k\in\{1,\dots,K\}$ and errors $E_1,\dots,E_K$ which are $N\times T$ random matrices such that\begin{equation}\label{model fac}
X_{kit} =  \sum_{j=1}^{r_N}\lambda_{ij}f_{tj}\delta_{kj}+E_{kit},\ k\in\{1,\dots,K\}.
\end{equation}
The role of the vectors $\delta, ..., \delta_K$ is to model the dependence between the regressors and the unobservables $\sum_{j=1}^{r_N}\lambda_{ij}f_{tj}\delta_{j}+E_{it}$. The structure that we impose can be seen as the generalisation to dimension $3$ (the third dimension being the one of variables) of the usual factor models for matrices as in \citet{bai2016econometric}. Such a modelling was already introduced in the psychometrics literature in \citet{carroll1970analysis} and \citet{harshman1970foundations}. The mathematical foundations behind this approach lie in the tensor decomposition literature, see \citet{kolda2009tensor} for a survey. 

We study a class of two-step estimators of the proposed model (\eqref{model 1} and \eqref{model fac}). In the first step  , the factors and the factor loadings are estimated. Then, in the second step the outcome is regressed on the covariates augmented by estimates of the factors and the factor loadings. We provide sufficient conditions on the first-step estimator under which the two-step estimator is asymptotically normal. We present assumptions under which a first-step estimator based on principal components analysis (henceforth PCA) satisfies these conditions. All the results are developed under an asymptotic regime where the sample size $N$ goes to infinity and $T$ is a function of $N$ going to infinity with $N$. Moreover, the number of factors is unknown and allowed to grow (possibly to infinity) with the sample size. Factors are not assumed to be strong. The proposed principal components augmented estimator exhibits better finite sample properties than alternatives in Monte-Carlo simulations.

When a strong factor assumption is imposed and the number of factors is assumed to be fixed, the proposed two-step estimator is found to be asymptotically normal under weaker conditions on $N$ and $T$ than for the factor-augmented estimator in \citet{greenaway2012asymptotic}. This suggests that augmenting the panel regression with estimates of the factor loadings leads to improved estimation properties. The estimator of Section 4.7.1 of \citet{beyhum2019square} is a special case of the two-step procedure of this paper. In this other article, a first-step estimator based on hard-thresholding of a nuclear-norm penalized estimator is used. The procedure is pivotal in the sense that it does not require knowledge of the variance of the error terms and that the thresholding level is data-driven.  The first-step estimator uses a penalty which level depend on the distribution of the operator norms of the errors while the approach with PCA that we develop here does not. It also relies on the fact that a compatibility constant is bounded away from $0$ with probability approaching $1$. Such an assumption is absent in the present paper. 

This paper is organized as follows.  The two-step estimator is introduced in Section \ref{sec:est}. Sufficient conditions for asymptotic normality are derived in Section \ref{sec:res}.  Section \ref{sec:PCA} is devoted to the analysis of the two-step procedure when PCA is used in the first step. Section \ref{sec:sim} describes our simulations.  All the proofs are deferred to the Appendix.

\noindent
\textbf{Preliminaries.} The transpose of a $N \times T$ matrix $A$ is written $A^{\top}$ and its trace is $\mathrm{tr}(A)$. Its $k^{th}$ singular value is $\sigma_k(A)$ and $\mathrm{rank}(A)$ is its rank. $A=\sum_{k=1}^{\mathrm{rank}(A)}\sigma_k(A)u_k(A)v_k(A)^\top$ is the singular value decomposition of $A$, where $\left\{u_k\left(A\right)\right\}_{k=1}^{\text{rank}\left(A\right)}$ is a family of orthonormal vectors of $\mathbb{R}^N$ and $\left\{v_k\left(A\right)\right\}_{k=1}^{\text{rank}\left(A\right)}$ is a family of orthonormal vectors of $\mathbb{R}^T$. The scalar product in the space of $N \times T$ matrices is $\left\langle A,B\right\rangle=\mathrm{tr}(A^{\top}B)$. The nuclear norm is $\left|A\right|_*=\sum_{k=1}^{\text{rank}(A)}\sigma_k(A)$, and the operator norm is $\left|A\right|_{\text{op}}=\sigma_1(A)=\max\limits_{h\in\mathbb{R}^T\ \text{s.t.}\ \left|h\right|_2=1} \left|Ah\right|_2$. For two integers, $N$ and $T$, $N \vee T$ is the maximum of $N$ and $T$ , $N \wedge T$ is the minimum of $N$ and $T$ and $\left \lfloor{N}\right \rfloor$ is the integer part of $N$. For $N\in \mathbb{N}$, $I_N$ is the identity matrix of size $N$. 

We consider sequences of data generating processes indexed by $N$. $T$ is a function of $N$ that goes to infinity with $N$. This paper studies an asymptotic where $N$ goes to infinity. For a probabilistic event $\mathcal{A}$, its complement is denoted $\mathcal{A}^c$ and we write that $\mathcal{A}$ happens with probability approaching $1$ or w.p.a $1$ if $\mathbb{P}(\mathcal{A})\to 1$.

\section{The estimator}
\label{sec:est}
The model can be rewritten in matrix form as $Y=\Pi_0+E_0$, $X_k=\Pi_k+E_k$ for $k\in\{1,\dots,K\}$, where $\Pi_{kit} = \sum_{j=1}^{r_N}\lambda_{ij}f_{tj}\delta_{kj}$ for $i\in \{1,\dots,N\}$ and $t\in \{1,\dots,T\}$, $\Pi_0= \sum_{k=1}^K\beta_k\Pi_k+\Gamma$, $\Gamma_{it} = \sum_{r=1}^{r_N}\lambda_{ir}f_{tr}\delta_{r}$ and $E_0= \sum_{k=1}^K\beta_kE_k+E$. Notice that $E_0$ and $E$ are different. $E$ is the remainder term in \eqref{model 1}, while $E_0$ is the remainder term in the expression of $Y$ as the sum of a term with a statistical factor structure and a remainder. Remark also that we do not assume that the error terms $E,E_0,\dots, E_K$ have mean zero, hence they can be the sum of an error term with mean zero and a small remainder as in \citet{beyhum2019square}.

Let $\Pi_u = (\Pi_0,\dots,\Pi_K)$,
$\Pi_v=((\Pi_0)^\top,\dots,(\Pi_K)^\top)$. For $z=u,v$, we denote by $P_z$ the projector on the vector space spanned by the columns of $\Pi_z$ and $M_z$ the projector on the orthogonal of the vector space spanned by the columns of $\Pi_z$. Let $r_z$ be the rank of $\Pi_z$. Note that $r_z\le r_N$, by definition. 

The proposed estimator is as follows. In a first step, one estimates $M_u$ and $M_v$ by estimators $\widehat{M}_u$ and $\widehat{M}_v$. From there, the estimator of $\beta$ is \begin{equation} \label{estimator}\widehat{\beta}\in \argmin{b\in \mathbb{R}^K}\left|\widehat{E}_0-\sum_{k=1}^{K}  b_{k}\widehat{E}_{k} \right|_2^2,\end{equation}
where $\widehat{E}_0=\widehat{M}_uY\widehat{M}_v$ and $\widehat{E}_k=\widehat{M}_uX_k\widehat{M}_v$ for $k\in\{1,\dots,K\}$. 

As argued in the introduction, the estimator \eqref{estimator} can be seen as the regression of the outcome on the regressors and estimated factor loadings and factors as shown in the following lemma. Let us introduce $\widehat{r}_u=\text{rank}\left(I_N-\widehat{M}_{u}\right)$, $\widehat{r}_v=\text{rank}\left(I_T-\widehat{M}_{v}\right)$ and $X_{it}=(X_{1it}, \dots, X_{Kit})^\top$.
\begin{lemma}\label{equivalence} Let $\{\widehat{\lambda_i}\}_{i=1}^{N}$ (resp. $\{\widehat{f_t}\}_{t=1}^{T}$) be a family of vectors in $\mathbb{R}^{\widehat{r}_u}$ (resp. $\mathbb{R}^{\widehat{r}_v}$) such that $\{(\widehat{\lambda}_{1j}\dots, \widehat{\lambda}_{Nj})^\top\}_{j=1}^{\widehat{r}_u}$ (resp. $\{(\widehat{f}_{1j}\dots, \widehat{f}_{Tj})^\top\}_{j=1}^{\widehat{r}_v}$) is a generating family of the orthogonal of the null space of $\widehat{M}_u$ (resp. $\widehat{M}_v$). Then, it holds that 
$$\widehat{\beta} \in \argmin{ b\in \mathbb{R}^K} \min\limits_{\begin{array}{c}
\phi_1,\dots, \phi_{T} \in \mathbb{R}^{\widehat{r}_u},\\
l_1,\dots, l_{N} \in \mathbb{R}^{\widehat{r}_v}\end{array}}
\sum_{i=1}^N \sum_{t=1}^T \left(Y_{it}-X_{it}^\top b- \widehat{\lambda}_i^\top \phi_t - l_i^\top \widehat{f}_t\right)^2.$$
\end{lemma}

\section{Sufficient assumptions for asymptotic normality}
\label{sec:res}
In this section, we present sufficient conditions for asymptotic normality of $\widehat{\beta}$ and consistent estimation of its asymptotic variance. The first assumption concerns the asymptotic behaviour of the error matrices. For a $N\times T$ matrix $A$, define $\widetilde{A}=M_uAM_v$.
\begin{assumption}
\label{E} The following holds:
\begin{enumerate}[\textup{(}i\textup{)}]  
 \item\label{Eii} There exists a $K\times K$ positive definite matrix $\Sigma$ such that, for $k,l\in \{1,\dots, K\}$, 
$\left\langle \widetilde{E}_k, \widetilde{E}_l\right\rangle/(NT) \xrightarrow{\P} \Sigma_{kl};$
 \item\label{Eiii} There exists $\sigma>0$ such that $\left|\widetilde{E}\right|_2^2/(NT)\xrightarrow{\P}\sigma^2$ and
$\left(\left\langle \widetilde{E}_k, \widetilde{E}\right\rangle\right)_{k=1}^{K}/\sqrt{NT} \xrightarrow{d} \mathcal{N}\left(0,\sigma^2\Sigma\right).$
 \end{enumerate}
 \end{assumption}
\noindent This assumption is similar to Assumption 9 (v) and (vi) in \citet{beyhum2019square}. The next lemma provides sufficient conditions for Assumption \ref{E}. 
\begin{lemma}
\label{lEG} Assume that
\begin{enumerate}[\textup{(}i\textup{)}]  
\item\label{lEGi} $\mathbb{E}\left[\left|P_u E\right|_2+ \left| EP_v\right|_2\right]+\sum_{k=1}^K\mathbb{E}\left[\left|P_u E_k\right|_2+ \left| E_kP_v\right|_2\right]=o_P(\sqrt{NT})$;
\item\label{lEGii}There exists a positive definite matrix $\Sigma$ such that, for $k,l\in\{1,\dots, K\}$, 
$\left\langle E_k,E_l\right\rangle/(NT)\xrightarrow{\P} \Sigma_{kl};$
\item\label{lEGiii} For $k\in\{1,\dots, K\}$, $\left\langle E_k, P_uE \right\rangle/\left|P_uE\right|_2=O_P\left(1\right)$ and $\left\langle E_k, M_uEP_v \right\rangle/\left|M_uEP_v\right|_2=O_P\left(1\right)$;
\item\label{lEGiv}There exists $\sigma >0$ such that $\left|E\right|_2^2/(NT)\xrightarrow{\P}\sigma^2$ and $\left(\left\langle E_k,E\right\rangle\right)_{k=1}^K/\sqrt{NT}\xrightarrow{d} \mathcal{N}\left(0,\sigma^2\Sigma\right)$.
 \end{enumerate}
Then, conditions \eqref{Eii} and \eqref{Eiii} in Assumption \ref{E} hold.
\end{lemma}\noindent
The next corollary gives an example of data generating process under which Assumption \ref{E} holds.
\begin{corollary}
\label{lE}
Let us assume that $E, E_1,\dots, E_k$ are independent, $r_u+r_v=o_P\left(\sqrt{N\wedge T}\right)$ and there exists $\sigma, \sigma_1,\dots,\sigma_k>0$ such that $\{E_{it}\}_{it}$ are i.i.d. $\mathcal{N}(0,\sigma^2)$ and $\{E_{kit}\}_{it}$ are i.i.d. $\mathcal{N}(0,\sigma_k^2)$. If also $(E, E_1,\dots, E_k)$ is independent of $(\Pi_u,\Pi_v)$, then Assumption \ref{E} holds.
\end{corollary}\noindent
The last set of conditions concerns the performance of the estimators of the projectors $\widehat{M}_u$ and $\widehat{M}_v$. Let $\{u_N\}_N$ and $\{v_N\}_N$ be real-valued sequences such that 
$\left|\widehat{M}_{u}-M_{u}\right|_2=O_P(u_N)$ and
$\left| \widehat{M}_{v}-M_{v}\right|_2=O_P(v_N)$. Let also $\{h_N\}_N$ and $\{\rho_N\}_N$ be real-valued sequences such that $\max\limits_{k\in\{0,\dots,K\}}\left|\Pi_k+E_k\right|_{2}= O_P(h_N)$ and $\max\limits_{k\in\{0,\dots,K\}}\left|E_k\right|_{\op}= O_P(\rho_N)$.
The estimators satisfy the following assumption.
\begin{assumption}
\label{Rates} The following holds:
\begin{enumerate}[\textup{(}i\textup{)}]  
\item \label{Ratesi} $\widehat{M}_{u}$ and $\widehat{M}_{v}$ are symmetric almost surely;
\item\label{Ratesiii} $\mathbb{P}(\widehat{r}_u = r_u) \to 1$ and $\mathbb{P}(\widehat{r}_v = r_v) \to 1$;
\item\label{Ratesiii.2} $u_N\vee v_N =o(1)$ and $h_N^2=O(NT)$;
\item\label{Ratesiv} $\sqrt{2r_N} (u_N\vee v_N)\rho_N^2 =o(NT)$;
\item\label{Ratesv} $u_Nv_Nh_N^2 =o(NT)$;
\item\label{Ratesvi} $\sqrt{2r_N} (u_N\vee v_N)\rho_N =o\left(\sqrt{NT}\right)$;
\item\label{Ratesvii} $u_Nv_N\rho_Nh_N =o\left(\sqrt{NT}\right)$;
\end{enumerate}
\end{assumption}
\noindent
This assumption plays a similar role as conditions (i) to (iv) in Assumption 9 in \citet{beyhum2019square}. It is difficult to understand the strength of Assumption \ref{Rates} without examples of $u_N$ and $v_N$ for specific first-step estimators. Hence, we discuss it in Section \ref{sec:PCA}, where we derive the properties of $\widehat{M}_u$ and $\widehat{M}_v$ when they are estimated by a method relying on PCA. The next theorem constitutes the main result of this paper.

\begin{theorem}[Asymptotic Normality] \label{AN}
Under assumptions \ref{E} and \ref{Rates}, we have 
$$\sqrt{NT}(\widehat{\beta}-\beta)\xrightarrow{d}\mathcal{N}\left(0,\sigma^2\Sigma^{-1}\right).$$
Also, for $k,l\in \{1,\dots,K\}$, $\widehat{\Sigma}_{kl}=\left\langle\widehat{E}_k,\widehat{E}_l\right\rangle/(NT)\xrightarrow{\P} \Sigma_{kl}$ and $\widehat{\sigma}^2=\left|\widehat{E}_{0}-\sum_{k=1}^{K}\widehat{\beta}_k\widehat{E}_{k}\right|_2^2/(NT)\xrightarrow{\P} \sigma^2.$ 
\end{theorem}
\noindent
\section{Estimation of the projectors using principal components analysis} \label{sec:PCA} 
\subsection{Strength of the factors} In this section, we discuss the estimation of the projectors using a method based on PCA. We make assumptions regarding the asymptotic behaviour of $\sigma_j(\Pi_z)$ for $z=u,v$. The purpose of this subsection is to show that there exists data generating processes (henceforth DGP) that generate various asymptotic behaviours of the singular values of $\Pi_z$. When $\sigma_j(\Pi_z)/\sqrt{NT}$ has a finite deterministic limit in probability, then we say that the $j^{th}$ factor is strong. The following lemma shows that there exists a wide variety of DGP under which such a strong factor assumption holds. Let $\Lambda=(\lambda_1,\dots,\lambda_N)$, $F=(f_1,\dots,f_N)$, and $\Delta = (\delta_0,\dots, \delta_K)$, where $\delta_0 =\delta+\sum_{k=1}^K\beta_k \delta_k$.
\begin{lemma}\label{spec}
Assume that $r_N$ is fixed and \begin{enumerate}[\textup{(}i\textup{)}]  
 \item\label{speci} There exists a $r_N\times r_N$ positive definite matrix $\Sigma_\Lambda$ such that $\Lambda\Lambda^\top/N \xrightarrow{\mathbb{P}}\Sigma_\Lambda$;
 \item\label{specii} There exists a $r_N\times r_N$ positive definite matrix $\Sigma_F$ such that $FF^\top/T \xrightarrow{\mathbb{P}}\Sigma_F$;
 \item\label{speciii} $\Delta\Delta^\top$ does not depend on $N$.
 \end{enumerate}
Then, for $z=u,v$, the ratio of the singular values of $\Pi_z$ and $\sqrt{NT}$ has a finite deterministic limit in probability.
\end{lemma}
If instead $\sigma_j(\Pi_z)/\alpha_{jN}$ has a finite deterministic limit for $\alpha_{jN}=o\left(\sqrt{NT}\right)$, then the $j^{th}$ factor is not strong. For a detailed discussion of the concept of non-strong factors, see \citet{pesaran2015time}. The following lemma shows how to generate non-strong factors and a growing number of factors in the case where $F$, $\Lambda$ and $\Delta$ are nonrandom.
\begin{lemma}\label{specweak}
Let  $\{\alpha_{jN}\}_N$ for $j\in \mathbb{N}$ be real-valued sequences with positive values. Maintain
 \begin{enumerate}[\textup{(}i\textup{)}]  
 \item\label{specweaki} $\Lambda$ is nonrandom and $(\Lambda\Lambda^\top)_{jj}=I_{r_N}$;
 \item\label{specweakii}  $F$ is nonrandom and $FF^\top$ is equal to the $r_N\times r_N$ diagonal matrix with coefficients $\alpha_{1r_N}^2, \dots,\alpha_{jN}^2$;
 \item\label{specweakiii} $\Delta\Delta^\top$ is a diagonal matrix such that $\alpha_{1 N}^2\left(\Delta\Delta^\top\right)_{11}\ge \dots\ge \alpha_{r_NN}^2\left(\Delta\Delta^\top\right)_{r_Nr_N}$.
 \end{enumerate}
Then, for $z=u,v$ and $r\in \mathbb{N}$, $\sigma_{j}(\Pi_z)=\alpha_{j N}\sqrt{\left(\Delta\Delta^\top\right)_{jj}}$.
\end{lemma}
Notice that the last two lemmas give sufficient conditions for Assumption 8 in \citet{beyhum2019square}.
\subsection{Convergence results} The econometrician can use different methods to estimate the projectors $M_u$ and $M_v$. The approach in \citet{beyhum2019square} relies on a nuclear-norm penalised estimator followed by hard-thresholding of the singular values. It has the advantage of being data-driven, in the sense that it does not use any knowledge of the number of factors or the variance of the errors. Another interesting and computationally advantageous procedure is the double IV estimator of \citet{gagliardini2017double}. In this paper, we focus the theoretical presentation on yet another method, based on the PCA. $r_u$ and $r_v$ are estimated via the eigenvalue ratio estimator from \citet{ahn2013eigenvalue}. For $z=u,v$, let us define
\begin{equation}
\widehat{r}_z\in \argmax_{j\in \left\{1,\dots,\left \lfloor{\sqrt{N\wedge T}}\right \rfloor \right\}} \frac{\sigma_{j}\left(Y_z\right)}{\sigma_{j+1}\left(Y_z\right)},\label{ERatio}
\end{equation}
where $Y_u=(Y,X_1,\dots, X_K)$ and $Y_v=(Y^\top,X_1^\top,\dots, X_K^\top)$. It may be that there exists $r\in \left\{1,\dots,\left \lfloor{\sqrt{N\wedge T}}\right \rfloor \right\}$, $\sigma_{j+1}\left(Y_z\right)=0$. To ensure that the estimators are defined, throughout this section, we use the convention that the division of a positive number by $0$ is equal to $\infty$.
The estimator in \citet{ahn2013eigenvalue} is of the form $\widehat{r}_z\in \argmax_{j\in \left\{1, \dots, \left \lfloor{d^*(N\wedge T)}\right \rfloor\right\}} \sigma_{j}\left(Y_z\right)/\sigma_{j+1}\left(Y_z\right)$, where $d^*\in (0,1]$. Therefore, the estimators in \eqref{ERatio} correspond to the one in \citet{ahn2013eigenvalue} for a particular choice of $d^*$. Our theoretical analysis is different from the one of \citet{ahn2013eigenvalue} because it allows for non-strong factors and a growing number of factors. Contrarily to the estimators in \citet{bai2002determining}, the advantage of the eigenvalue ratio estimator is that it does not require to choose a penalty level. To ensure consistency of the eigenvalue ratio estimator, we make the following assumption. Let $E_u=(E_0,\dots, E_K)$ and $E_v=(E_0^\top,\dots, E_K^\top)$.
\begin{assumption}[Eigenvalue Ratio] \label{ER}For $z=u,v$, it holds that $r_z\le \sqrt{N\wedge T}$ almost surely, $\left|E_z \right|_\op=o_P\left(\sigma_{r_z}\left(\Pi_z\right)\right)$ and there exists $C<1$ such that
\begin{equation}\notag \mathbb{P}\left(\left(\max_{j\in\left\{1,\dots, r_z-1\right\}} \frac{\sigma_j\left(\Pi_z\right)}{\sigma_{j+1}\left(\Pi_z\right)}\right)\vee \left(\max_{j \in \left\{r_z+1,\dots, \left \lfloor{\sqrt{N\wedge T}}\right \rfloor\right\}}\frac{ \left\vert E_z \right|_{\op}}{\sigma_{r_z+j}\left(E_z\right)}\right)\le C \frac{\sigma_{r_z}\left(\Pi_z\right)}{\sigma_{2r_z+1}\left(E_z\right)} \right)\to 1.\end{equation}
\end{assumption}

\noindent 
Let us give sufficient conditions fo Assumption \ref{ER}.
\begin{lemma} \label{exER} For $z=u,v$, assume that $r_z\le \sqrt{N\wedge T}$, $\left|E_z \right|_\op=O_P\left(\sqrt{N \vee T}\right)$, $\left|E_z\right|_2^2/(NT)$ has a finite deterministic limit in probability and there exists a sequence $\{z_N\}_N$ such that $\sigma_{r_z}(\Pi_z)=O_P(z_N)$, $\sqrt{N\vee T}=o\left(z_N\right)$ and $\max_{j\in\left\{1,\dots, r_z-1\right\}} \sigma_j\left(\Pi_z\right)/\sigma_{j+1}\left(\Pi_z\right)=o_P\left(z_N/\sqrt{N\vee T}\right)$,  then Assumption \ref{ER} holds.
\end{lemma}
\noindent This Lemma shows that our assumption allows for non-strong factors and a growing number of factors. The condition $\max_{j\in\left\{1,\dots, r_z-1\right\}} \sigma_j\left(\Pi_z\right)/\sigma_{j+1}\left(\Pi_z\right)=o_P\left(z_N/\sqrt{N\vee T}\right)$ implies that the singular values of $\Pi_z$ cannot decrease too quickly with $j\in\{1,\dots, r_z\}$. The assumption that $\left|E_z \right|_\op=O_P\left(\sqrt{N \vee T}\right)$ is standard in the panel data literature and holds under flexible cross-sectional and serial correlations. For a detailed discussion, see Appendix A.1 in \citet{moon2015linear}. Let us now state the main result regarding the eigenvalue ratio estimator. 
\begin{lemma}\label{LERatio} Under Assumption \ref{ER}, we have $\mathbb{P}\left(\widehat{r}_u=r_u,\widehat{r}_v=r_v\right)\to 1$.\end{lemma}

Given the estimators $\widehat{r}_u$ and $\widehat{r}_v$, we set
$\widehat{M}_u=I_N-\sum_{j=1}^{\widehat{r}_u}u_j\left(Y_u\right)u_j\left(Y_u\right)^{\top}$ and 
$\widehat{M}_v=I_T-\sum_{j=1}^{\widehat{r}_v}u_j\left(Y_v\right)u_j\left(Y_v\right)^{\top}$.
\noindent Then, we have the following theorem which states the rates of convergence of the estimators of the projectors.
\begin{theorem} \label{Proj}For $z=u,v$, if $\mathbb{P}\left(\widehat{r}_z=r_z\right)\to 1$, we have
$\left|\widehat{M}_z-M_z\right|_{2} =O_P\left(\sqrt{r_z}\left|E_z\right|_{\op}/\sigma_{r_z}\left(\Pi_z\right)\right)$.
\end{theorem}
\subsection{Examples}\label{subsec:ex}

Let us now show how Assumption \ref{Rates} can hold under different assumptions on the singular values of $\Pi_u,\Pi_v, E_u$ and $E_v$. In both examples, we assume that, for $z=u,v$, $\left|E_z \right|_\op=O_P\left(\sqrt{N \vee T}\right)$ and $\left|E_z\right|_2^2/(NT)$ has a deterministic finite limit. The assumption on the errors implies that we can choose $\rho_N= \sqrt{N\vee T}$ because $\left|E_k\right|_{\op}\le \left|E_u\right|_{\op}$.

\noindent
\textbf{Example 1.} In this first example, we assume that $r_N$ is fixed and that the strong factor assumption holds, that is, for $z=u,v$ and $j\in\{1,\dots, r\}$, $\sigma_j(\Pi_z)/\sqrt{NT}$ has a finite deterministic limit. This implies that we can choose $h_N=\sqrt{NT}$ because $\left|\Pi_u\right|_2=O_P\left(\sqrt{NT}\right)$ and that Assumption \ref{ER} holds by Lemma \ref{exER}. Theorem \ref{Proj} yields $u_N=v_N=1/\sqrt{N\wedge T} $. All conditions in Assumption \ref{Rates} except \eqref{Ratesvii} are satisfied whatever the value of $N$ and $T$. Condition \eqref{Ratesvii} holds if $\sqrt{N\vee T}/(N\wedge T)=o(1)$. The latter correponds to the condition for asymptotic normality of the debiased estimator in \citet{bai2009panel} and is weaker than the conditions for asymptotic normality in  \citet{greenaway2012asymptotic}.

\noindent
\textbf{Example 2.} In this case, we assume that $r_u=r_v=r_N$ can grow with the sample size, and that, for $z=u,v$ and $j\in\{1,\dots, r_N\}$, $\sqrt{r_N}\sigma_j(\Pi_z)/\sqrt{NT}$ has a finite deterministic limit. This is a case with non-strong factors and a growing number of factors. This implies that we can choose $h_N=\sqrt{NT}$ because $\left|\Pi_u\right|_2=O_P\left(\sqrt{NT}\right)$ and that Assumption \ref{ER} holds by Lemma \ref{exER}.  From Theorem \ref{Proj}, we obtain $u_N=v_N=r_N/\sqrt{N\wedge T}$. Conditions \eqref{Ratesi}-\eqref{Ratesiii.2} hold for any value of $N$, $T$ and $r_N$. For \eqref{Ratesiv}- \eqref{Ratesvi} to hold, it is enough that $r_N^{\frac32}/(N\wedge T)=o(1)$. Finally, condition \eqref{Ratesvii} is satisfied if $r_N^2\sqrt{N\vee T}/(N\wedge T)=o(1)$.

\section{Simulations}
\label{sec:sim}

We consider a data generating process with a single regressor and two factors:
\begin{align*}
Y_{it}&= X_{1it}+\lambda_{i1}f_{t1}+\lambda_{i2}f_{t2}+E_{it},\\
X_{1it}&=\frac{1}{2}\lambda_{i1}f_{t1}+\lambda_{i2}f_{t2}+E_{1it},
\end{align*} where $f_{tl}$, $\lambda_{il}$, $E_{1it}$, and $E_{it}$ for all indices are mutually independent, $f_{tl}\sim\mathcal{N}(1/2,1)$, $\lambda_{il}\sim\mathcal{N}(1,1)$ and $E_{1it}, \dots, E_{Kit}$ and $E_{it}$ are standard normals. The matrix $X_1$ has a statistical factor structure with a low-rank component of rank 2. Recall that $\widehat{\beta}^{LS}\in \argmin{b\in \mathbb{R}}\left|Y-bX_1\right|_2^2$ is the least-squares estimator of the linear regression of the outcome on the regressors. $\widehat{\beta}^{FA}\in \argmin{b\in \mathbb{R}}\left|Y\widehat{M}_v-bX_1\widehat{M}_v\right|_2^2$ is the factor augmented regression estimator where $\widehat{M}_v$ is computed as in Section \ref{sec:PCA}. $\widetilde{\beta}^{(1)}$ and $\widetilde{\beta}^{(2)}$ are the two-stage estimators of Section 4.7.1 in \citet{beyhum2019square}. They are computed as in the simulations of that paper, without using within transforms. $\widetilde{\beta}^{(2)}$ uses Bai's estimator as a second stage while $\widetilde{\beta}_2$ uses the approach of this paper with two projectors and a first-step based on hard-thresholding of a nuclear-norm penalized estimator. Finally, $\widehat{\beta}^{PCA}$ is the estimator \eqref{estimator}, using the procedure of Section \ref{sec:PCA} as the first-stage.

Tables \ref{fig:Cov} and \ref{fig:Cov1} compare the performance of the estimators in terms of mean squared error (henceforth MSE), bias, standard error (henceforth std) and coverage of 95\% confidence intervals, for different sample sizes. The coverage is not reported for $\widehat{\beta}_{LS}$ because the latter is not asymptotically normal for the DGP that we consider. We use 7300 Monte-Carlo replications which allows for an accuracy of $\pm0.005$ with 95\% for the coverage probabilities of 95\% confidence intervals. In this simulation exercise, our estimator exhibits better finite samples properties than the studied alternatives.

\begin{table}[htbp]
  \centering
  \caption{$N=T=50$}
\label{fig:Cov} 
{\small
             \begin{tabular}{|c|c|c|c|c|c|}
  \hline
 & $\widehat{\beta}^{LS}$  & $\widehat{\beta}^{FA}$& $\widetilde{\beta}^{(1)}$ & $\widetilde{\beta}^{(2)}$ & $\widehat{\beta}^{PCA}$\\
  \hline
   MSE  & 0.884 & 0.004 & 0.14&  0.13 & 0.004 \\
   bias  & 0.939 &-0.011 &0.321&  0.275& 0.012\\
   std  &   0.055 & 0.191 & 0.023 &   0.234&  0.063 \\
      coverage  &    & 0.75& 0.22 &   0.37&  0.90\\
  \hline
      \end{tabular}  }
\end{table}
%\begin{table}[!ht]

\begin{table}[htbp]
  \centering
  \caption{$N=T=150$}
    \label{fig:Cov1}
    {\small
             \begin{tabular}{|c|c|c|c|c|c|}
  \hline
 & $\widehat{\beta}^{LS}$  & $\widehat{\beta}^{FA}$& $\widetilde{\beta}^{(1)}$ & $\widetilde{\beta}^{(2)}$ & $\widehat{\beta}^{PCA}$ \\
  \hline
   MSE  & 0.887 & $10^{-4}$ & 4 $10^{-5}$&  4 $10^{-5}$ & 4 $10^{-5}$ \\
   bias  & 0.9414&-0.007& -5 $10^{-5}$&  -8 $10^{-6}$ & -5 $10^{-5}$\\
   std  &   0.031& 0.007& 0.007 &   0.007& 0.007 \\
         coverage  &  & 0.79  & 0.95 &   0.95&  0.95 \\
  \hline
        \end{tabular}      
        }
\end{table}

\bibliographystyle{plainnat}
\bibliography{refPCA}
\section*{Appendix}
\subsection*{Proof of Lemma \ref{equivalence}.}
Let $b\in \mathbb{R}^k$, $Y_i=(Y_{i1},\dots, Y_{iT})^\top$ and $X_i=(X_{i1},\dots, X_{iT})^\top$ for $i\in\{1,\dots, N\}$, $Y_t=(Y_{1t},\dots, Y_{Nt})^\top$ and $X_t=(X_{1t},\dots, X_{Nt})^\top$ for $t\in\{1,\dots, T\}$ and 
$$\varphi(b)=\min\limits_{\begin{array}{c}
\phi_1,\dots, \phi_{T} \in \mathbb{R}^{\widehat{r}_u},\\
l_1,\dots, l_{N} \in \mathbb{R}^{\widehat{r}_v}\end{array}}
\sum_{i=1}^N \sum_{t=1}^T \left(Y_{it}-X_{it}^\top b- \widehat{\lambda}_i^\top \phi_t - l_i^\top \widehat{f}_t\right)^2.$$
By algebra, we have $$\varphi(b) =\min\limits_{\begin{array}{c}
\phi_1,\dots, \phi_{T} \in \mathbb{R}^{\widehat{r}_u},\\
l_1,\dots, l_{N} \in \mathbb{R}^{\widehat{r}_v}\end{array}}
\sum_{i=1}^N  \left|Y_{i}-X_{i}b-  \left(\phi_1,\dots,\phi_T\right)^\top\widehat{\lambda}_i- \left(\widehat{f}_1,\dots, \widehat{f}_T\right)^\top l_i\right|^2_2.$$
Then, by definition of $\widehat{M}_v$, it holds
$$\varphi(b)=\min\limits_{\begin{array}{c}
\phi_1,\dots, \phi_{T} \in \mathbb{R}^{\widehat{r}_u}\end{array}}
\sum_{i=1}^N  \left|\widehat{M}_v\left(Y_{i}-X_{i} b-  \left(\phi_1,\dots,\phi_T\right)^\top\widehat{\lambda}_i\right)\right|^2_2.$$
Because $\widehat{M}_v$ is symmetric, this implies
\begin{equation}\label{avantdernierM}
\varphi(b) =  \min\limits_{\begin{array}{c}
\phi_1,\dots, \phi_{T} \in \mathbb{R}^{\widehat{r}_u}\end{array}}
 \left|\left(Y-\sum_{k=1}^K b_kX_{k}-\left(\widehat{\lambda}_1,\dots,\widehat{\lambda}_N\right)^\top  \left(\phi_1,\dots,\phi_T\right)\right)\widehat{M}_v\right|^2_2.\end{equation}
Next, by definition of $\widehat{M_u}$, we obtain 
$$\varphi(b) \ge   
 \left|\widehat{M}_u\left(Y-\sum_{k=1}^K b_kX_{k}\right)\widehat{M}_v\right|^2_2.$$
 Hence, because the value of $$ \min\limits_{\begin{array}{c}
\phi_1,\dots, \phi_{T} \in \mathbb{R}^{\widehat{r}_u}\end{array}}
 \left|\left(Y-\sum_{k=1}^K b_kX_{k}-\left(\widehat{\lambda}_1,\dots,\widehat{\lambda}_N\right)^\top  \left(\phi_1,\dots,\phi_T\right)\right)\widehat{M}_v\right|^2_2$$ is $\left|\widehat{M}_u\left(Y-\sum_{k=1}^K b_kX_{k}\right)\widehat{M}_v\right|^2_2$ when $\left(\widehat{\lambda}_1,\dots,\widehat{\lambda}_N\right)^\top \phi_t =  \left(I_N-\widehat{M}_u\right)(Y_t-X_t b)$, we get 
 
 $$\varphi(b) =   
 \left|\widehat{M}_u\left(Y-\sum_{k=1}^K b_kX_{k}\right)\widehat{M}_v\right|^2_2.$$

\subsection*{Proof of Lemma \ref{lEG}. }
\noindent \textbf{Proof that Assumption \ref{E} \eqref{Eii} holds.}
 For $k\in \{1,\dots,K\}$, we have
$$E_k-M_{u}E_kM_{v}=
P_{u}E_k+M_{u}E_kP_{v}. $$
By Markov's inequality and the fact that $M_u$ is a projector, we have $$\left|P_{u}E_k\right|_2+\left|M_{u}E_kP_v\right|_2\le \left|P_{u}E_k\right|_2+\left|E_kP_v\right|_2=O_P\left(\mathbb{E}\left[\left|P_{u}E_k\right|_2+\left|E_kP_v\right|_2\right]\right).$$
By condition \eqref{lEGi} in Lemma \ref{lEG}, this yields $\left|P_{u}E_k\right|_2+\left|M_{u}E_kP_v\right|_2=o_P\left(\sqrt{NT}\right)$. 
This implies $\left\vert M_{u}E_k M_{v}-E_k\right\vert_2=o_P\left(\sqrt{NT}\right)$. 
By the Cauchy-Schwarz inequality, we obtain $$\left\langle M_{u}E_k M_{v},E_l\right\rangle-\left\langle E_k,E_l\right\rangle=\left\langle M_{u}E_kM_{v}-E_k,E_l\right\rangle=o_P\left(NT\right)$$
because $\left| E_k \right|_2=O_P\left(\sqrt{NT}\right)$. We get $$
\frac{1}{NT}\left\langle\widetilde{E}_k,\widetilde{E}_l\right\rangle= \frac{1}{NT}\left\langle M_{u}E_k M_{v},E_l\right\rangle= \frac{1}{NT}\left\langle E_k ,E_l\right\rangle +o_P(1)
 \xrightarrow{\P} \Sigma_{kl},$$
by condition \eqref{lEGii} in Lemma \ref{lEG}.\\
\noindent
 \textbf{Proof that Assumption \ref{E} \eqref{Eiii} holds.} The proof of $\left|\widetilde{E}\right|_2^2/(NT)\xrightarrow{\P}\sigma^2$ is similar to the proof that Assumption \ref{E} \eqref{Eii} holds.
By conditions \eqref{lEGi} and \eqref{lEGiii} in Lemma \ref{lEG}, we have $$\left\langle P_{u}E_k,E\right\rangle=O_P(1)\left| P_{u}E_k\right|_2= O_P(1)o_P\left(\sqrt{NT}\right)=o_P\left(\sqrt{NT}\right)$$ and  similarly $\left\langle M_{u}E_kP_v,E\right\rangle=o_P\left(\sqrt{NT}\right).$ Next, this yields
$$
\left\langle E_k ,E\right\rangle-\left\langle M_{u}E_kM_{v} ,E\right\rangle\
=\left\langle P_{u}E_k,E\right\rangle+\left\langle M_{u}E_kP_{v} ,E\right\rangle=o_P\left(\sqrt{NT}\right).
$$
We obtain that $$
\frac{1}{\sqrt{NT}}\left(\left\langle\widetilde{E}_k,\widetilde{E}\right\rangle\right)_{k=1}^K= \frac{1}{\sqrt{NT}}\left(\left\langle M_{u}E_kM_{v} ,E\right\rangle\right)_{k=1}^K
=\frac{1}{\sqrt{NT}}\left(\left\langle E_k ,E\right\rangle\right)_{k=1}^K+o_P(1) \xrightarrow{d} \mathcal{N}\left(0,\sigma^2\Sigma\right),$$
by condition \eqref{lEGiv} in Lemma \ref{lEG}.

\subsection*{Proof of Corollary \ref{lE}.}
Let us prove that the assumptions of Proposition \ref{lEG} are satisfied. \eqref{lEGii} and \eqref{lEGiv} in Proposition \ref{lEG} are direct consequences of the weak law of large numbers and the central limit theorem. Concerning \eqref{lEGi} in Proposition \ref{lEG},  for $k\in \{1,\dots, K\}$, by Lemma A.3 in \citet{giraud2014introduction} and the fact that $E,E_1,\dots,E_K$ are independent of $\Pi_u$ and $\Pi_v$, we have $$\mathbb{E}\left[\left\vert P_{u}E_k\right\vert_{2}^2 \right]=\sum_{t=1}^T\mathbb{E}\left[\sum_{i=1}^N\left( P_{u}E_k\right)_{it}^2\right]=\sum_{t=1}^T\mathbb{E}\left[\mathbb{E}\left[\left.\sum_{i=1}^N\left( P_{u}E_k\right)_{it}^2\right\vert P_u\right]\right]=Tr_u \sigma_k^2=o\left(\sqrt{NT}\right)$$
Similarly, one can show that $\mathbb{E}\left[\left\vert E_kP_v\right\vert_{2}^2 \right]=o\left(\sqrt{NT}\right)$. In the same manner, we obtain that $\mathbb{E}\left[\left|P_u E\right|_2+ \left| EP_v\right|_2\right]=o(\sqrt{NT})$. To prove \eqref{lEGiii}, just notice that conditionally on $P_uE$, we have 
$\left\langle P_uE, E_k,  \right\rangle/\left|P_uE\right|_2\sim\mathcal{N}(0,\sigma_k^2)$, hence,
for any $M\ge 0$, it holds that
$$\mathbb{P}\left(\left.\frac{\left|\left\langle P_{u}E_k ,E\right\rangle\right|}{\left\vert P_{u}E\right\vert_2\sigma_k}> M\right\vert  P_uE \right)\le 2(1-\Phi^{-1}(M)),$$
where $\Phi$ is the cumulative distribution function of a $\mathcal{N}(0,1)$ distribution. This implies 
$$\mathbb{P}\left(\frac{\left|\left\langle P_{u}E_k ,E\right\rangle\right|}{\left\vert P_{u}E\right\vert_{2}\sigma_k}>M\right)\le 2(1-\Phi^{-1}(M)).$$
Therefore, we obtain that $\left\langle P_{u}E_k ,E\right\rangle/\left\vert P_{u}E\right\vert_{2}=O_P(1)$. The proof that $\left\langle E_k ,M_uEP_v\right\rangle/\left\vert M_{u}EP_v\right\vert_{2}=O_P(1)$ is the same.

\subsection*{Proof of Theorem \ref{AN}.} \noindent \textbf{Proof of asymptotic normality.}

\noindent Because $\widehat{M}_u$ and $\widehat{M}_v$ are symmetric, a solution to \eqref{estimator} satisfies, for $l=1,\dots,K$, 
$$\left\langle \widehat{M}_uX_l\widehat{M}_v ,Y-\sum_{k=1}^K\widehat{\beta}_kX_k\right\rangle=0,$$ hence
\begin{align*}
&\left\langle M_{u}X_lM_{v} ,+E+\sum_{k=1}^K\left(\beta_k-\widehat{\beta}_k\right)X_k \right\rangle\\
&=\left\langle\left(M_u-\widehat{M}_u\right)X_lM_v ,E+\sum_{k=1}^K\left(\beta_k-\widehat{\beta}_k\right)X_k\right\rangle\\
&\quad+\left\langle M_uX_l\left(M_v-\widehat{M}_v\right) ,E+\sum_{k=1}^K\left(\beta_k-\widehat{\beta}_k\right) X_k\right\rangle\\
&\quad-\left\langle \left(M_u-\widehat{M}_u\right) X_l\left(M_v-\widehat{M}_v\right) ,\Gamma+E+\sum_{k=1}^K\left(\beta_k-\widehat{\beta}_k\right) X_k \right\rangle, 
\end{align*}
so
\begin{align}
&\sum_{k=1}^K\left(\beta_k-\widehat{\beta}_k\right)\Bigg(
\left\langle M_uX_lM_v ,X_k\right\rangle
-\left\langle\left(M_u-\widehat{M}_u\right)X_lM_v ,X_k \right\rangle\notag\\
& \quad 
-\left\langle M_uX_l\left(M_v-\widehat{M}_v\right) ,X_k \right\rangle\notag\\
& \quad +\left\langle \left(M_u-\widehat{M}_u\right) X_l\left(M_v-\widehat{M}_v\right) ,X_k \right\rangle\Bigg)+\left\langle\left(M_u-\widehat{M}_u\right)X_lM_v ,+E\right\rangle\notag\\
&\quad +\left\langle M_uX_l\left(M_v-\widehat{M}_v\right) ,+E\right\rangle\notag\\
&\quad -\left\langle \left(M_u-\widehat{M}_u\right) X_l\left(M_v-\widehat{M}_v\right),\Gamma+E\right\rangle.\label{emain}
\end{align}
Let us show that $\left\langle M_uX_lM_v ,X_k\right\rangle$, which by Assumption \ref{E} \eqref{Eii} diverges like $NT$, is the high-order term %in the bracket 
multiplying $\left(\beta_k-\widehat{\beta}_k\right)$ in \eqref{emain}. This also yields the consistency of the estimator of the %variance-
covariance matrix.  For a matrix $M$ and $r\in\N$, let us define $|M|_{2,r}^2=\sum_{k=1}^{r}\sigma_k(M)^2$.
By symmetry of the projectors, Theorem C.5 in \cite{giraud2014introduction}, and Assumption \ref{Rates} \eqref{Ratesiii} $\left(\text{which implies $\rank\left(M_u-\widehat{M}_u\right)\le 2r_N$ w.p.a. $1$}\right)$, we have
\begin{align*}
&\left|\left\langle\left(M_u-\widehat{M}_u\right)X_lM_{v} ,X_k \right\rangle\right|\\
&\le
\left|M_u-\widehat{M}_u\right|_2\left|X_l M_{v}X_k^\top\right|_{2,2r_N}\\
%\left|\left(M_u-\widehat{M}_u\right)X_lM_{v}\right|_*\left|X_kM_{v\}\right|_{\op}\\ 
%&\le \left(\sqrt{2r_N}+o_P(1)\right)\left|\left(M_u-\widehat{M}_u\right)X_lM_v\right|_2\left|X_kM_v\right|_{\op}\\
&\le \left(\sqrt{2r_N}+o_P(1)\right)\left|M_u-\widehat{M}_u\right|_2\left|X_lM_v\right|_{\op}\left|X_kM_v\right|_{\op}\\
&
=O_P\left(\sqrt{2r_N} u_N\rho_N^2\right)= o_P(NT) \quad (\text{by Assumption \ref{Rates} \eqref{Ratesiv}}).
\end{align*}
We bound similarly $\left|\left\langle M_uX_l\left(M_v-\widehat{M}_v\right) ,X_k\right\rangle\right|$,
and, for the fourth term, use that
%, using that $\rank\left(M_u-\widehat{M}_u\right)\le 2r_N$ w.p.a. 1 and Lemma C2 item 1 in \cite{giraud2014introduction},
\begin{align}
\notag &\left|\left\langle \left(M_u-\widehat{M}_u\right) X_l\left(M_v-\widehat{M}_v\right) ,X_k \right\rangle
\right|\\
\notag &\le \left| \left(M_u-\widehat{M}_u\right) X_l\left(M_v-\widehat{M}_v\right)\right|_2 \left|X_k \right|_{2}\\
\label{double02}&=O_P\left(u_Nv_Nh_N^2\right)= o_P(NT) \quad (\text{by Assumption \ref{Rates} \eqref{Ratesv}}).
\end{align}
Let us consider now the quantities on the right-hand side in \eqref{emain}. Notice that because $E =E_0-\sum_{k=1}^K\beta_kE_k$, it holds that $\left|E\right|_{\op}=O_P(\rho_N)$. Proceeding like above,
we have
\begin{align*}
&\left|\left\langle\left(M_u-\widehat{M}_u\right)X_lM_v ,E \right\rangle\right|\notag\\
&\le
\left|M_u-\widehat{M}_u\right|_2\left|X_l M_vE^\top\right|_{2,2r_N}\notag\\
&\le 
\left(\sqrt{2r_N}+o_P(1)\right)u_N\left|E_l\right|_{\op}\left|E\right|_{\op}\\
& =O_P\left(\sqrt{2r_N} u_N\rho_N^2\right)=o_P(\sqrt{NT}) \quad (\text{by Assumption \ref{Rates} \eqref{Ratesvi}}).
\end{align*}
and treat similarly $\left\langle M_uX_l\left(M_v-\widehat{M}_v\right) ,E\right\rangle$.
With the same arguments as in \eqref{double02}, the absolute value of the last term of \eqref{emain} is smaller than
$u_Nv_N\left|X_l\right|_{2}\left|\Gamma+E\right|_{2}$, which is an $O_P(u_Nv_N\rho_Nh_N)=o_P\left(\sqrt{NT}\right)$ because $\Gamma +E= Y-\sum_{k=1}^K\beta_kX_k $.

Let us now look at the first terms on the left-hand side and on the right-hand side of \eqref{emain}. By Assumption \ref{Rates} \eqref{Ratesvi},
for all $k,l\in\{1,\dots,K\}$, we have $\left\langle M_uX_lM_v ,X_k\right\rangle=\left\langle M_uE_lM_v ,E_k\right\rangle+o_P(NT)$. Hence
because of Assumption \ref{E} \eqref{Eii}, $\left\langle M_uX_lM_v ,X_k\right\rangle$ are the high-order terms on the left-hand side of \eqref{emain}. Similarly, by Assumption \ref{E} \eqref{Eiii}, the high-order terms on the right-hand side of \eqref{emain} are $\left\langle M_uE_lM_v ,E\right\rangle$. 
As a result, $\widehat{\beta}$ is asymptotically equivalent to the ideal estimator $\overline{\beta}$
\begin{equation}\label{ewanideal}
\overline{\beta}\in\argmin{\beta\in\mathbb{R}^K}\left|M_u\left(Y-\sum_{k=1}^K\beta_kX_k\right)M_v\right|_2^2.
\end{equation}
Hence, we obtain by usual arguments that $\sqrt{NT}(\widehat{\beta}-\beta)\xrightarrow{d}\mathcal{N}\left(0,\sigma\Sigma^{-1}\right)$. \\

\noindent \textbf{Proof of the consistency of $\widehat{\sigma}$.} We use
\begin{align*} 
NT\widehat{\sigma}^2&= \left\langle Y-\sum_{k=1}^{K} \widehat{\beta}_kX_k ,\widehat{M}_u\left(Y-\sum_{k=1}^{K} \widehat{\beta}_kX_k\right)\widehat{M}_v\right\rangle \\
&=  \left\langle Y-\sum_{k=1}^{K} \widehat{\beta}_kX_k , \left(\widehat{M}_u-M_u  \right)\left(Y-\sum_{k=1}^{K} \widehat{\beta}_kX_k\right)\left(\widehat{M}_v-M_v \right)\right\rangle\\
&\quad + \left\langle Y-\sum_{k=1}^{K} \widehat{\beta}_kX_k , \left(\widehat{M}_u-M_u  \right)\left(Y-\sum_{k=1}^{K} \widehat{\beta}_kX_k\right)M_v \right\rangle\\
&\quad +  \left\langle Y-\sum_{k=1}^{K} \widehat{\beta}_kX_k ,M_u\left(Y-\sum_{k=1}^{K} \widehat{\beta}_kX_k\right)\left(\widehat{M}_v-M_v \right)\right\rangle\\
&\quad +  \left\langle Y-\sum_{k=1}^{K} \widehat{\beta}_kX_k ,M_u \left(Y-\sum_{k=1}^{K} \widehat{\beta}_kX_k\right)M_v \right\rangle.\\
\end{align*}
Now, by the Cauchy-Schwarz inequality,
\begin{align*}
& \left|\left\langle Y-\sum_{k=1}^{K} \widehat{\beta}_kX_k , \left(\widehat{M}_u-M_u  \right)\left(Y-\sum_{k=1}^{K} \widehat{\beta}_kX_k\right)\left(\widehat{M}_v-M_v \right)\right\rangle \right|\\
 &\le 
\left|Y-\sum_{k=1}^{K} \widehat{\beta}_kX_k \right|_2\left|\left(\widehat{M}_u-M_u  \right)\left(Y-\sum_{k=1}^{K} \widehat{\beta}_kX_k\right)\left(\widehat{M}_v-M_v \right) \right|_2\\
&\le \left|Y-\sum_{k=1}^{K} \widehat{\beta}_kX_k\right|_2^2\left|\widehat{M}_u-M_u \right|_2 \left|\widehat{M}_v-M_v \right|_2\\
&=O_P(h_N^2 u_Nv_N)\quad \text{(by the fact that $\widehat{\beta}-\beta=o_P(1)$)}\\
&=o_P(NT)\quad \text{(by Assumption \ref{Rates} \eqref{Ratesiii.2})}).
\end{align*}
Similarly, we can show that $$\left\langle Y-\sum_{k=1}^{K} \widehat{\beta}_kX_k , \left(\widehat{M}_u-M_u  \right)\left(Y-\sum_{k=1}^{K} \widehat{\beta}_kX_k\right)M_v \right\rangle=o_P(NT)$$ and $$ \left\langle Y-\sum_{k=1}^{K} \widehat{\beta}_kX_k,M_u\left(Y-\sum_{k=1}^{K} \widehat{\beta}_kX_k\right)\left(\widehat{M}_v-M_v \right)\right\rangle=o_P(NT).$$
Hence, we have \begin{align*}&NT\widehat{\sigma}^2\\
&= \left\langle Y-\sum_{k=1}^{K} \widehat{\beta}_kX_k ,M_u \left(Y-\sum_{k=1}^{K} \widehat{\beta}_kX_k\right)M_v \right\rangle+o_P(NT)\\
&=  \left\langle Y-\sum_{k=1}^{K} \left(\widehat{\beta}_k-\beta_k\right)X_k -\sum_{k=1}^{K} \beta_kX_k ,M_u \left(Y-\sum_{k=1}^{K}  \left(\widehat{\beta}_k-\beta_k\right)X_k-\sum_{k=1}^{K+1} \beta_kX_k\right)M_v \right\rangle+o_P(NT)\\
&= \left\langle \sum_{k=1}^{K} \left(\widehat{\beta}_k-\beta_k\right)X_k  ,M_u \left(\sum_{k=1}^{K}  \left(\widehat{\beta}_k-\beta_k\right)X_k\right)M_v \right\rangle+ \left\langle E_0  ,M_u \left(\sum_{k=1}^{K}  \left(\widehat{\beta}_k-\beta_k\right)X_k\right)M_v \right\rangle\\
&\quad +\left\langle \sum_{k=1}^{K} \left(\widehat{\beta}_k-\beta_k\right)X_k  ,M_u  E_{K} M_v \right\rangle+ \left|\widetilde{E}\right|_2^2+o_P(NT).
\end{align*}
Now, by the Cauchy-Schwarz inequality, Assumption \ref{Rates} and the fact that $\widehat{\beta}-\beta=o_P(1)$, one can show that
\begin{align*}\left\langle \sum_{k=1}^{K} \left(\widehat{\beta}_k-\beta_k\right)X_k  ,M_u \left(\sum_{k=1}^{K}  \left(\widehat{\beta}_k-\beta_k\right)X_k\right)M_v \right\rangle=o_P(NT)\\
 \left\langle E  ,M_u \left(\sum_{k=1}^{K}  \left(\widehat{\beta}_k-\beta_k\right)X_k\right)M_v \right\rangle=o_P(NT);\\
\left\langle \sum_{k=1}^{K} \left(\widehat{\beta}_k-\beta_k\right)X_k  ,M_u  E M_v \right\rangle=o_P(NT).
\end{align*}
We conclude the proof using Assumption \ref{E}.
\subsection*{Proof of Lemma \ref{spec}.}
Let $\Lambda =(\lambda_1,\dots\lambda_N)$. For $t\in\{1,\dots ,T\}$ and $k\in\{0,\dots, K\}$, we use the notation $\psi_{tk}=\left(f_{t1}\delta_{k1},\dots,f_{tr_N}\delta_{kr_N}\right)^\top$. We also introduce 
$\Psi=(\psi_{10},\dots,\psi_{T0},\dots, \psi_{1K},\dots,\psi_{TK}).$ It holds that, for $j,j'\in \{1,\dots, r_N\}$, $\left(\Psi\Psi^\top\right)_{jj'}/T=(\Delta\Delta^\top)_{jj'}\left(FF^\top\right)_{jj'}/T$.
Therefore,  $\Lambda\Lambda^\top \Psi\Psi^\top/(NT)$ converges in probability to $\Sigma_{\Lambda}\Sigma_{\Delta F}$, where, for $j,j'\in \{1,\dots, r_N\}$, $\left(\Sigma_{\Delta F}\right)_{jj'} =(\Delta\Delta^\top)_{jj'} \left(\Sigma_{F}\right)_{jj'}$.
Next, let $U=(u_1\left(\Pi_u\right)),\dots,u_{r_N}\left(\Pi_u\right))$, $V=(v_1\left(\Pi_u\right),\dots, v_{r_N}\left(\Pi_u\right))$ and $D$ be the $r_N\times r_N$ diagonal matrix for which $D_{jj}=\sigma_j\left(\Pi_u\right)$. We have $UDV^\top=\Lambda^\top\Psi$, which implies $UD^2U^\top=\Lambda^\top\Psi\Psi^\top\Lambda$. This yields $\Lambda UD^2= \Lambda\Lambda^\top \Psi\Psi^\top\Lambda U$.
On the event $\mathcal{E}=\{\mathrm{rank}(\Lambda U)=r_N\}$, we obtain $ \Lambda\Lambda^\top \Psi\Psi^\top=\Lambda UD^2 (\Lambda U)^{-1}$. Therefore, the diagonal elements of $D^2/(NT)$ are the eigenvalues of $\Lambda\Lambda^\top \Psi\Psi^\top/(NT)$ on $\mathcal{E}$. Because $\Lambda\Lambda^{\top}/N$ converges in probability to a positive definite matrix, the set of full rank matrices is an open set and the determinant is a continuous mapping, we have $\mathbb{P}(\mathrm{rank}(\Lambda)=r_N)\to 1$, which  implies $\mathbb{P}(\mathcal{E})\to 1$. For $j\in\{1,\dots, r_N\}$ and $\xi >0$, we get
\begin{align*}
&\mathbb{P}\left(\left|\sigma_j\left(\frac{D^2}{NT}\right)-\sigma_j\left(\Sigma_{\Lambda}\Sigma_{\Delta F}\right)\right|\le \xi \right)\\
& \ge \mathbb{P}\left(\left\{\left|\sigma_j\left(\frac{D^2}{NT}\right)-\sigma_j\left(\Sigma_{\Lambda}\Sigma_{\Delta F}\right)\right|\le \xi\right\}\cap \mathcal{E} \right)\\
&=\mathbb{P}\left(\left\{\left|\sigma_j\left(\Lambda\Lambda^\top \Psi\Psi^\top/(NT)\right)-\sigma_j\left(\Sigma_{\Lambda}\Sigma_{\Delta F}\right)\right|\le \xi\right\}\cap \mathcal{E} \right)\to 1,
\end{align*}
where the last statement holds because $\Lambda\Lambda^\top \Psi\Psi^\top/(NT)\xrightarrow{\mathbb{P}}\Sigma_{\Lambda}\Sigma_{\Delta F}$, $A\in \mathbb{R}^{r_N\times r_N}\mapsto (\sigma_1(A),\dots, \sigma_{r_N}(A)).$ is a continuous mapping and $\mathbb{P}(\mathcal{E})\to 1$.
\subsection*{Proof of Lemma \ref{specweak}.}
We only prove the result for $\Pi_u$, the proof for $\Pi_v$ being similar. We use the same notations as in the proof of Lemma \ref{spec}. We have, for $j,j '\in \{1,\dots, r_N\}$, $\left(\Psi\Psi^\top\right)_{jj'}=(FF^\top)_{jj'}(\Delta\Delta^\top)_{jj'}=0$ if $j\ne j'$ and $\left(\Psi\Psi^\top\right)_{jj}=(FF^\top)_{jj}(\Delta\Delta^\top)_{jj} = \alpha_{jN}^2(\Delta\Delta^\top)_{jj}$ if $j=j'$. Therefore, $\Lambda\Lambda^\top \Psi\Psi^\top$ is the diagonal matrix with diagonal coefficients $\alpha_{1N}^2(\Delta\Delta^\top)_{11},\dots, \alpha_{r_NN}^2(\Delta\Delta^\top)_{r_Nr_N}$. Because $\Lambda$ has full rank, $ \Lambda\Lambda^\top \Psi\Psi^\top=\Lambda UD^2 (\Lambda U)^{-1}$ and, therefore, the diagonal coefficients of $D^2$ are $\alpha_{1N}^2(\Delta\Delta^\top)_{11},\dots, \alpha_{r_NN}^2(\Delta\Delta^\top)_{r_Nr_N}$.

\subsection*{Results on PCA}
Let us consider a $N \times T$ random matrix $A$.  We do not observe $A$ but 
$\widetilde{A}=A+Z,$
where $Z$ is an $N \times T$ random matrix. Let $r$ be the rank of $A$. $A=\sum_{j=1}^r \sigma_ju_jv_j^\top $ is the singular value decomposition of $A$, where $\sigma_1\ge \dots\ge \sigma_r\ge 0$ and $\left\{u_1,\dots,u_r\right\}$ and $\left\{v_1,\dots,v_r\right\}$ are orthonormal families of $\R^N$ and $\R^T$, respectively. With similar notations, $\widetilde{A}=\sum_{j=1}^{\widetilde{r} }\widetilde{\sigma}_j \widetilde{u}_j\widetilde{v}_j^\top $ is the singular value decomposition of $\widetilde{A}$ and $\widetilde{r}$ is the rank of $\widetilde{A}$. $Z=\sum_{j=1}^{N \wedge T}\sigma_j(Z) u_j(Z)v_j(Z)^\top $ is a singular value decomposition of $Z$. $T=T(N)$ is a function of $N$ going to $\infty$ when $N\to \infty$ and and the asymptotic setting is such that $N\to \infty$. For $s\in\{1,\dots, N\wedge T\}$, we consider the following estimators of $A$ and $P$, $\widehat{A}_{s}=\sum_{j=1}^{s}\widetilde{\sigma}_j\widetilde{u}_j\widetilde{v}_j^\top $ and $\widehat{P}_s=\sum_{j=1}^{s} \widetilde{u}_j\widetilde{u}_j^{\top}$. Let also $\widehat{M}_{s}=I_N-\widehat{P}_{s}$.
\begin{lemma}\label{l1}
$\left\vert \widehat{A}_r-A\right\vert_{\op}
\le 2\left\vert Z \right\vert_{\rm op}.$
\end{lemma}
\textbf{Proof.}
We have
$
\left\vert \widehat{A}_r-A\right\vert_{\op}
=\left\vert \widehat{A}_r-\widetilde{A}+\widetilde{A}-A\right\vert_{\op}\le \left\vert \sum_{j=r+1}^{N\wedge T}\widetilde{\sigma}_j \widetilde{u}_j\widetilde{v}_j^\top \right\vert_{\rm op}+\left\vert Z \right\vert_{\rm op} =\widetilde{\sigma}_{r+1}+\left\vert Z\right\vert_{\rm op}.
$
Now, by Weyl's inequality (Theorem C.6 in \citet{giraud2014introduction}), it holds that $\widetilde{\sigma}_{r+1} \le \left|\widetilde{A} -A\right|_{\op}=\left\vert Z \right\vert_{\rm op}$.
\hfill $\Box$
\begin{lemma} \label{l2}
We have 
$\left| \widehat{P}_r -P\right \vert_2\le 4\sqrt{2r}\frac{\left|Z\right|_{\op}}{\sigma_r}$ almost surely.
\end{lemma}
\textbf{Proof.} Following the proof of Proposition 10 in \citet{beyhum2019square}, we obtain $\left| \widehat{P}_r -P\right \vert_2^2\le 2\left|\widehat{M}_rA\right|_{2}^2/\sigma_r^2.$ We conclude using  $$\left \vert \widehat{M}_rA \right\vert_2= \left \vert \widehat{M}_r\left(\widehat{A}_r-A\right) \right \vert_2 \le \left \vert\widehat{A}_r-A \right \vert_{2}\le \sqrt{2r} \left \vert\widehat{A}_r-A \right \vert_{\op}\le \sqrt{2r} 2\left| Z\right|_{\op},$$
by Lemma \ref{l1} and the fact that $\widehat{M}_r$ is a projector. \hfill $\Box$

\begin{lemma} The following holds:
\label{l4}
\begin{enumerate}[\textup{(}i\textup{)}]  
\item\label{l4i} For $j\in\{1,\dots, r\}$, $\sigma_j- \left|\widehat{A}_r-A\right|_{\op}\le \widetilde{\sigma}_j\le \sigma_j+ \left|\widehat{A}_r-A\right|_{\op}$;
\item\label{l4ii}For $j\in \{r+1, N\wedge T - r\}$,  $\sigma_{r+j}\left(Z\right)\le \widetilde{\sigma}_j\le \left|Z \right|_{\op}$.\end{enumerate}
\end{lemma}
\textbf{Proof.}
\eqref{l4i} follows from the fact that $\left| \widetilde{\sigma}_j-\sigma_j \right|\le \left|\widehat{A}_r-A\right|_{\op}$ by Weyl's inequality. Weyl's inequality also yields $\widetilde{\sigma}_j\le \left|\widetilde{A}-A\right|_{\op}=\left|Z \right|_{\op}$, which implies the right-hand side of \eqref{l4ii}. To show the left-hand side of \eqref{l4ii}, from (7.3.13) in \citet{horn2012matrix}, we obtain $\sigma_{r+j}(Z)\le \widetilde{\sigma}_{j-1}+\sigma_{r+1}= \widetilde{\sigma}_j$.
\hfill $\Box$
\begin{lemma}\label{suff} Let $Z$ be a $N \times T$ random matrix and $r\in\left\{1,\dots,\left\lfloor{\sqrt{N\wedge T}}\right \rfloor\right\}$. Assume that $\left|Z \right|_\op=O_P\left(\sqrt{N \vee T}\right)$ and there exists $v> 0$ such that $\left|Z\right|_2^2/(NT)\xrightarrow{\P} v^2$. Then, we have $\sigma_{2\left \lfloor{\sqrt{N\wedge T}}\right \rfloor}(Z)>0$ w.p.a. $1$ and $ \max\limits_{j\in\left\{1,\dots, \left \lfloor{\sqrt{N\wedge T}}\right \rfloor\right\}} \left\vert Z\right|_{\op}/\sigma_{r+j}(Z)=O_P(1).$\end{lemma}
\textbf{Proof.}
We have 
$$\frac{\left|Z \right|_2^2}{NT}\le \frac{2 \left \lfloor{\sqrt{N\wedge T}}\right \rfloor}{NT}\left|Z\right|_{\op}^2+ \frac{N \wedge T}{NT} \sigma_{2\left \lfloor{\sqrt{N\wedge T}}\right \rfloor}(Z)^2.$$
Thus, we obtain 
$$\frac{\left|Z \right|_2^2}{NT}-\frac{2\left \lfloor{\sqrt{N\wedge T}}\right \rfloor}{NT}\left|Z\right|_{\op}^2\le  \frac{N \wedge T}{NT} \sigma_{2\left \lfloor{\sqrt{N\wedge T}}\right \rfloor}(Z)^2.$$
Using $\left|Z\right|_2^2/(NT)\xrightarrow{\P} v^2$ and $\left|Z \right|_\op=O_P\left(\sqrt{N \vee T}\right)$, we get $$\frac{\left|Z \right|_2^2}{NT}-\frac{2\left \lfloor{\sqrt{N\wedge T}}\right \rfloor}{NT}\left|Z\right|_{\op}^2\xrightarrow{\mathbb{P}}v^2$$ and, therefore,
$$\P\left(\frac{\sigma_{2\left \lfloor{\sqrt{N\wedge T}}\right \rfloor}(Z)}{\sqrt{N \vee T}} \ge \frac{v}{2}\right)\to 1.$$ Hence, we have
$$ \P\left(\frac{\left|Z \right|_{\op}}{\sigma_{2\left \lfloor{\sqrt{N\wedge T}}\right \rfloor}(Z)} \le \frac{2\left|Z \right|_{\op}}{\sqrt{N \vee T}v}\right)\to 1.$$
Therefore, we obtain
$$ \frac{\left|Z\right|_{\op}}{\sigma_{2\left \lfloor{\sqrt{N\wedge T}}\right \rfloor}(Z)} =O_P\left(\frac{\left|Z \right|_{\op}}{\sqrt{N \vee T}}\right)=O_P(1).$$
This leads to 
$$ \max_{j\in \{1, \dots,\left \lfloor{\sqrt{N\wedge T}}\right \rfloor\}}\frac{ \left\vert Z\right|_{\op}}{\sigma_{r+j}(Z)}\le \frac{\left|Z\right|_{\op}}{\sigma_{2\left \lfloor{\sqrt{N\wedge T}}\right \rfloor}(Z)}=O_P(1) .$$
\hfill $\Box$

\subsection*{Proof of Lemma \ref{exER}.} Because $\sqrt{N\vee T}=o(z_N)$, it holds that $\left|E_z\right|_{\op}=O_P(\sigma_{r_z}(\Pi_z))$. Then, we have $\sigma_{2r_z+1}(E_z)\le \left|E_z\right|_{\op} =O_P\left(\sqrt{N\vee T}\right)$ which implies $z_N/\sqrt{N\vee T}=O_P\left(\sigma_{r_z}(\Pi_z)/\sigma_{2r_z+1}(E_z)\right)$. Moreover, by Lemma \ref{suff}, we have $ \max\limits_{j\in\left\{1,\dots, \left \lfloor{\sqrt{N\wedge T}}\right \rfloor\right\}}\left\vert E_z\right|_{\op}/\sigma_{r_z+j}(E_z)=O_P(1)$. Because $$\max_{j\in\left\{1,\dots, r_z-1\right\}} \sigma_r\left(\Pi_z\right)/\sigma_{r+1}\left(\Pi_z\right)=o_P\left(z_N/\sqrt{N\vee T}\right),$$
we obtain
$$\mathbb{P}\left(\left(\max_{j\in\left\{1,\dots, r_z-1\right\}} \frac{\sigma_r\left(\Pi_z\right)}{\sigma_{r+1}\left(\Pi_z\right)}\right)\vee \left(\max_{j \in \left\{r_z+1,\dots, \left \lfloor{\sqrt{N\wedge T}}\right \rfloor\right\}}\frac{ \left\vert E_z \right|_{\op}}{\sigma_{r_z+j}\left(E_z\right)}\right)\le C \frac{\sigma_{r_z}\left(\Pi_z\right)}{\sigma_{2r_z+1}\left(E_z\right)} \right)\to 1.
$$
\subsection*{Proof of Lemma \ref{LERatio}.}
To prove Lemma \ref{LERatio}, let us show that $\mathbb{P}\left(\max_{j\in\{1,\dots, r_z-1\}}\frac{\sigma_j(Y_z)}{\sigma_{j+1}(Y_z)}<\frac{\sigma_r(Y_z)}{\sigma_{r+1}(Y_z)}\right)\to 1$.
Take $j\in\{1,\dots, r_z-1\}$, by Lemma \ref{l4} \eqref{l4i} and Lemma \ref{l1}, we have $\sigma_j(\Pi_z)- 2\left|E_z\right|_{\op}\le \sigma_j(Y_z)\le \sigma_j(\Pi_z)+2 \left|E_z\right|_{\op}$ . Then, on the event $\mathcal{A}=\left\{\sigma_{r_z}(\Pi_z)>2\left|E_z\right|_{\op}\right\}$, we obtain
\begin{equation} \label{endessous}\frac{\sigma_j(Y_z)}{\sigma_{j+1}(Y_z)}\le \frac{\sigma_j(\Pi_z)+2\left|E_z\right|_{\op}}{\sigma_{j+1}(\Pi_z)-2\left|E_z\right|_{\op}}=  \frac{\sigma_j(\Pi_z)}{\sigma_{j+1}(\Pi_z)} \frac{1+\frac{2\left|E_z\right|_{\op}}{\sigma_{j}(\Pi_z)}}{1-\frac{2\left|E_z\right|_{\op}}{\sigma_{j+1}(\Pi_z)}}\le  \frac{\sigma_j(\Pi_z)}{\sigma_{j+1}(\Pi_z)} \frac{1+\frac{2\left|E_z\right|_{\op}}{\sigma_{r_z}(\Pi_z)}}{1-\frac{2\left|E_z\right|_{\op}}{\sigma_{r_z}(\Pi_z)}},
\end{equation}
where the last equality is because $\sigma_j(\Pi_z)\ge \sigma_{j+1}(\Pi_z)\ge \sigma_{r_z}(\Pi_z)$.
Also, by Lemma \ref{l4}, on $\mathcal{A}$, it holds that
\begin{equation}\label{milieu}
\frac{\sigma_{r_z}(Y_z)}{\sigma_{r_z+1}(Y_z)}\ge \frac{\sigma_{r_z}(\Pi_z)-2\left|E_z\right|_{\op}}{\sigma_{2r_z+1}(E_z)}
= \frac{\sigma_{r_z}(\Pi_z)}{\sigma_{2r_z+1}(E_z)} \left(1-\frac{2\left|E_z\right|_{\op}}{\sigma_{r_z}(\Pi_z)}\right).
\end{equation}
Let us call $\mathcal{B}$ the event 
$$\left\{\left(1-\frac{2\left|E_z\right|_{\op}}{\sigma_{r_z}(\Pi_z)}\right)^2/\left(1+\frac{2\left|E_z\right|_{\op}}{\sigma_{r_z}(\Pi_z)}\right)> C\right\},$$
where $C$ is the constant in Assumption \ref{ER}.
We have
\begin{align*}
&\P\left(\max_{j\in\{1,\dots, r_z-1\}}\frac{\sigma_j(Y_z)}{\sigma_{j+1}(Y_z)}<\frac{\sigma_{r_z}(Y_z)}{\sigma_{r_z+1}(Y_z)}\right) \\
&\ge \P\left(\left\{\max_{j\in\{1,\dots, r_z-1\}}\frac{\sigma_j(Y_z)}{\sigma_{j+1}(Y_z)}<\frac{\sigma_{r_z}(Y_z)}{\sigma_{r_z+1}(Y_z)}\right\}\cap  \mathcal{A}\cap \mathcal{B}\right)\\
& \ge \P\left(\left\{\max_{j\in\{1,\dots, r_z-1\}}\frac{\sigma_j(\Pi_z)}{\sigma_{j+1}(\Pi_z)}<\frac{\left(1-\frac{2\left|E_z\right|_{\op}}{\sigma_{r_z}}(\Pi_z)\right)^2}{1+\frac{2\left|E_z\right|_{\op}}{\sigma_{r_z}(\Pi_z)}}\frac{\sigma_{r_z}(\Pi_z)}{\sigma_{2r_z+1}(E_z)}\right\}\cap \mathcal{A}\cap \mathcal{B}\right)\quad\text{(by \eqref{endessous} and \eqref{milieu})}\\
& \ge \P\left(\left\{\max_{j\in\{1,\dots, r_z-1\}}\frac{\sigma_j(\Pi_z)}{\sigma_{j+1}(\Pi_z)}<C\frac{\sigma_{r_z}(\Pi_z)}{\sigma_{2r_z+1}(E_z)}\right\}\cap \mathcal{A}\cap \mathcal{B}\right)\to 1,
\end{align*}
where the last statement holds because $\mathbb{P}(\mathcal{A})\to 1$, $\mathbb{P}(\mathcal{B})\to 1$ (given that $\left|E_z\right|_{\op} =O_P(\sigma_{r_z}(\Pi_z))$) and $$\P\left(\max_{j\in\{1,\dots, r_z-1\}}\frac{\sigma_j(\Pi_z)}{\sigma_{j+1}(\Pi_z)}<C\frac{\sigma_{r_z}(\Pi_z)}{\sigma_{2r_z+1}(E_z)}\right)\to 1$$ by Assumption \ref{ER}.
Next, let us show that, $\P\left(\max_{j\in \{r_z+1, \dots,\left \lfloor{\sqrt{N\wedge T}}\right \rfloor\}}\frac{\sigma_j(Y_z)}{\sigma_{j+1}(Y_z)}<\frac{\sigma_{r_z}(Y_z)}{\sigma_{r_z+1}(Y_z)}\right)\to 1$. By Lemma \ref{l4} \eqref{l4ii}, we have, for all $j>r_z$,
\begin{equation} \label{audessus} \frac{\sigma_j(Y_z)}{\sigma_{j+1}(Y_z)}\le \frac{ \left\vert E_z \right|_{\op}}{\sigma_{r_z+j}(E_z)}.\end{equation} 
Let $\mathcal{C}= \left\{ 1-\frac{2\left|E_z\right|_{\op}}{\sigma_{r_z}(\Pi_z)}>C\right\}$. This implies that 
\begin{align*}
&\P\left(\max_{j\in\left\{r_z+1, \dots, \left\lfloor{\sqrt{N \wedge T}}\right\rfloor\right\}}\frac{\sigma_j(Y_z)}{\sigma_{j+1}(Y_z)}<\frac{\sigma_{r_z}(Y_z)}{\sigma_{r_z+1}(Y_z)}\right) \\
&\ge \P\left(\left\{\max_{j\in \left\{r_z+1, \dots, \left\lfloor{\sqrt{N \wedge T}}\right\rfloor\right\}}\frac{\sigma_j(Y_z)}{\sigma_{j+1}(Y_z)}<\frac{\sigma_{r_z}(Y_z)}{\sigma_{r_z+1}(Y_z)}\right\}\cap \mathcal{A}\cap \mathcal{C}\right)\\
& \ge \P\left(\left\{\frac{ \left\vert E_z \right|_{\op}}{\sigma_{r_z+j}(E_z)}< \left(1-\frac{2\left|E_z\right|_{\op}}{\sigma_{r_z}(\Pi_z)}\right)\frac{\sigma_{r_z}(\Pi_z)}{\sigma_{2r_z+1}(E_z)}\right\}\cap \mathcal{A}\cap \mathcal{C}\right) \quad\text{(by \eqref{milieu} and \eqref{audessus})}\\
& \ge \P\left(\left\{\frac{ \left\vert E_z \right|_{\op}}{\sigma_{r+j}(E_z)}<C\frac{\sigma_{r_z}(\Pi_z)}{\sigma_{2r_z+1}(E_z)}\right\}\cap  \mathcal{A}\cap \mathcal{C}\right)\to 1,
\end{align*}
where the last statement holds because $\mathbb{P}(\mathcal{A})\to 1$, $\mathbb{P}(\mathcal{C})\to 1$ (given that $\left|E_z\right|_{\op} =O_P(\sigma_{r_z}))$) and $$\P\left(\frac{ \left\vert E_z \right|_{\op}}{\sigma_{r_z+j}(E_z)}<C\frac{\sigma_{r_z}(\Pi_z)}{\sigma_{r_z+1}(\Pi_z)}\right)\to 1$$ by Assumption \ref{ER}.
In the end, we obtain $$\P\left(\left(\max\limits_{j\in\{1,\dots, r_z-1\}}\frac{\sigma_j(Y_z)}{\sigma_{j+1}(Y_z)}\right)\vee\left(\max\limits_{j\in \left\{r_z+1, \dots, \left\lfloor{\sqrt{N \wedge T}}\right\rfloor\right\}}\frac{\sigma_j(Y_z)}{\sigma_{j+1}(Y_z)}\right)<\frac{\sigma_{r_z}(Y_z)}{\sigma_{2r_z+1}(Y_z)}\right)\to 1,$$ which concludes the proof.

\subsection*{Proof of Theorem \ref{Proj}.} We denote $\mathcal{A}= \left\{ \widehat{r}_z=r_z\right\}$. We have 
$$\P\left( \left| \widehat{M}_{z} -M_z\right \vert_2 \le 4 \sqrt{2r_z} \frac{\left|E_z\right|_{\op}}{\sigma_{r_z}(\Pi_z)}\right)
\ge  \P\left(\left\{ \left| \widehat{P}_z -P_z\right \vert_2\le  4\sqrt{2r_z} \frac{\left| E_z\right|_{\op}}{\sigma_{r_z}(\Pi_z)}\right\} \cap \mathcal{A}\right)=\P(\mathcal{A})
\to 1,
$$by Lemma \ref{l2}.

\end{document}